\magnification=1200

\hsize=11.25cm    
\vsize=18cm       
\parindent=12pt   \parskip=5pt     

\hoffset=.5cm   
\voffset=.8cm   

\pretolerance=500 \tolerance=1000  \brokenpenalty=5000

\catcode`\@=11

\font\eightrm=cmr8         \font\eighti=cmmi8
\font\eightsy=cmsy8        \font\eightbf=cmbx8
\font\eighttt=cmtt8        \font\eightit=cmti8
\font\eightsl=cmsl8        \font\sixrm=cmr6
\font\sixi=cmmi6           \font\sixsy=cmsy6
\font\sixbf=cmbx6

\font\tengoth=eufm10 
\font\eightgoth=eufm8  
\font\sevengoth=eufm7      
\font\sixgoth=eufm6        \font\fivegoth=eufm5

\skewchar\eighti='177 \skewchar\sixi='177
\skewchar\eightsy='60 \skewchar\sixsy='60

\newfam\gothfam           \newfam\bboardfam

\def\tenpoint{
  \textfont0=\tenrm \scriptfont0=\sevenrm \scriptscriptfont0=\fiverm
  \def\rm{\fam\z@\tenrm}
  \textfont1=\teni  \scriptfont1=\seveni  \scriptscriptfont1=\fivei
  \def\oldstyle{\fam\@ne\teni}\let\old=\oldstyle
  \textfont2=\tensy \scriptfont2=\sevensy \scriptscriptfont2=\fivesy
  \textfont\gothfam=\tengoth \scriptfont\gothfam=\sevengoth
  \scriptscriptfont\gothfam=\fivegoth
  \def\goth{\fam\gothfam\tengoth}
  
  \textfont\itfam=\tenit
  \def\it{\fam\itfam\tenit}
  \textfont\slfam=\tensl
  \def\sl{\fam\slfam\tensl}
  \textfont\bffam=\tenbf \scriptfont\bffam=\sevenbf
  \scriptscriptfont\bffam=\fivebf
  \def\bf{\fam\bffam\tenbf}
  \textfont\ttfam=\tentt
  \def\tt{\fam\ttfam\tentt}
  \abovedisplayskip=12pt plus 3pt minus 9pt
  \belowdisplayskip=\abovedisplayskip
  \abovedisplayshortskip=0pt plus 3pt
  \belowdisplayshortskip=4pt plus 3pt 
  \smallskipamount=3pt plus 1pt minus 1pt
  \medskipamount=6pt plus 2pt minus 2pt
  \bigskipamount=12pt plus 4pt minus 4pt
  \normalbaselineskip=12pt
  \setbox\strutbox=\hbox{\vrule height8.5pt depth3.5pt width0pt}
  \let\bigf@nt=\tenrm       \let\smallf@nt=\sevenrm
  \normalbaselines\rm}

\def\eightpoint{
  \textfont0=\eightrm \scriptfont0=\sixrm \scriptscriptfont0=\fiverm
  \def\rm{\fam\z@\eightrm}
  \textfont1=\eighti  \scriptfont1=\sixi  \scriptscriptfont1=\fivei
  \def\oldstyle{\fam\@ne\eighti}\let\old=\oldstyle
  \textfont2=\eightsy \scriptfont2=\sixsy \scriptscriptfont2=\fivesy
  \textfont\gothfam=\eightgoth \scriptfont\gothfam=\sixgoth
  \scriptscriptfont\gothfam=\fivegoth
  \def\goth{\fam\gothfam\eightgoth}
  
  \textfont\itfam=\eightit
  \def\it{\fam\itfam\eightit}
  \textfont\slfam=\eightsl
  \def\sl{\fam\slfam\eightsl}
  \textfont\bffam=\eightbf \scriptfont\bffam=\sixbf
  \scriptscriptfont\bffam=\fivebf
  \def\bf{\fam\bffam\eightbf}
  \textfont\ttfam=\eighttt
  \def\tt{\fam\ttfam\eighttt}
  \abovedisplayskip=9pt plus 3pt minus 9pt
  \belowdisplayskip=\abovedisplayskip
  \abovedisplayshortskip=0pt plus 3pt
  \belowdisplayshortskip=3pt plus 3pt 
  \smallskipamount=2pt plus 1pt minus 1pt
  \medskipamount=4pt plus 2pt minus 1pt
  \bigskipamount=9pt plus 3pt minus 3pt
  \normalbaselineskip=9pt
  \setbox\strutbox=\hbox{\vrule height7pt depth2pt width0pt}
  \let\bigf@nt=\eightrm     \let\smallf@nt=\sixrm
  \normalbaselines\rm}

\tenpoint

\def\pc#1{\bigf@nt#1\smallf@nt}         \def\pd#1 {{\pc#1} }

\catcode`\;=\active
\def;{\relax\ifhmode\ifdim\lastskip>\z@\unskip\fi
\kern\fontdimen2  -1.2 \fontdimen3 \string;}

\catcode`\:=\active
\def:{\relax\ifhmode\ifdim\lastskip>\z@\unskip\fi\penalty\@M\ \fi\string:}

\catcode`\!=\active
\def!{\relax\ifhmode\ifdim\lastskip>\z@
\unskip\fi\kern\fontdimen2  -1.1 \fontdimen3 \string!}

\catcode`\?=\active
\def?{\relax\ifhmode\ifdim\lastskip>\z@
\unskip\fi\kern\fontdimen2  -1.1 \fontdimen3 \string?}

\frenchspacing

\def\raggedbottom{\topskip 10pt plus 36pt\r@ggedbottomtrue}

\def\pointir{\unskip . --- \ignorespaces}

\def\Medbreak{\vskip-\lastskip\medbreak}

\long\def\th#1 #2\enonce#3\endth{
   \Medbreak\noindent
   {\pc#1} {#2\unskip}\pointir{\it #3}\smallskip}

\def\decale#1{\smallbreak\hskip 28pt\llap{#1}\kern 5pt}
\def\decaledecale#1{\smallbreak\hskip 34pt\llap{#1}\kern 5pt}
\def\puce{\smallbreak\hskip 6pt{$\scriptstyle\bullet$}\kern 5pt}

\def\eqalign#1{\null\,\vcenter{\openup\jot\m@th\ialign{
\strut\hfil$\displaystyle{##}$&$\displaystyle{{}##}$\hfil
&&\quad\strut\hfil$\displaystyle{##}$&$\displaystyle{{}##}$\hfil
\crcr#1\crcr}}\,}

\catcode`\@=12

\showboxbreadth=-1  \showboxdepth=-1

\newcount\numerodesection \numerodesection=1
\def\section#1{\bigbreak
 {\bf\number\numerodesection.\ \ #1}\nobreak\medskip
 \advance\numerodesection by1}

\mathcode`A="7041 \mathcode`B="7042 \mathcode`C="7043 \mathcode`D="7044
\mathcode`E="7045 \mathcode`F="7046 \mathcode`G="7047 \mathcode`H="7048
\mathcode`I="7049 \mathcode`J="704A \mathcode`K="704B \mathcode`L="704C
\mathcode`M="704D \mathcode`N="704E \mathcode`O="704F \mathcode`P="7050
\mathcode`Q="7051 \mathcode`R="7052 \mathcode`S="7053 \mathcode`T="7054
\mathcode`U="7055 \mathcode`V="7056 \mathcode`W="7057 \mathcode`X="7058
\mathcode`Y="7059 \mathcode`Z="705A


\def\hfl#1#2#3{\smash{\mathop{\hbox to#3{\rightarrowfill}}\limits
^{\textstyle#1}_{\textstyle#2}}}

\def\ogoth{{\goth o}}

\def\pgoth{{\goth p}}

\def\Q{{\bf Q}}
\def\Qp{\Q_p}

\def\Z{{\bf Z}}
\def\Zp{\Z_p}
\def\F{{\bf F}}
\def\Fp{{\F_{\!p}}}

\def\Gal{\mathop{\rm Gal}\nolimits}
\def\Ker{\mathop{\rm Ker}\nolimits}

\def\to{\rightarrow}

\def\normressym(#1,#2)_#3{\displaystyle\left({#1,#2\over#3}\right)}

\def\mod{\mathop{\rm mod.}\nolimits}
\def\pmod#1{\;(\mod#1)}

\newcount\refno 
\long\def\ref#1:#2<#3>{                                        
\global\advance\refno by1\par\noindent                              
\llap{[{\bf\number\refno}]\ }{#1} \pointir{\it #2} #3\goodbreak }

\def\citer#1(#2){[{\bf\number#1}\if#2\empty\relax\else,\ {#2}\fi]}

\def\\#1{\hbox{$\oldstyle#1$}}

\newbox\bibbox
\setbox\bibbox\vbox{\bigbreak
\centerline{{\pc BIBLIOGRAPHIC} {\pc REFERENCES}}

\ref {\pc BOREVICH} (Z.) \& {\pc SHAFAREVICH} (I.):
Number theory,
<Academic Press, New York-London, 1966, x+435 pp.>
\newcount\bosha \global\bosha=\refno

\ref{\pc DALAWAT} (C.):
Local discriminants, kummerian extensions, and elliptic curves,
<Journal of the Ramanujan Mathematical Society, {\bf 25} (2010) 1,
pp.~25--80. Cf.~arXiv\string:0711.3878v1.>      
\newcount\locdisc \global\locdisc=\refno

\ref{\pc HILBERT} (D.):
Die Theorie der algebraischen Zahlk{\"o}rper,
<Jahresbericht der Deutschen Mathematikervereinigung, {\bf 4} (1897),
pp.~175--546. $=$
Ges.\ Abhandlungen I, pp.~63--363.>
\newcount\zahlbericht \global\zahlbericht=\refno

\ref{\pc KATO} (K.), {\pc KUROKAWA} (N.) \& {\pc SAITO} (T.):
Number theory 1, Fermat's dream,
<American Mathematical Society, Providence, 2000, xvi+154 pp.>
\newcount\kks \global\kks=\refno

\ref{\pc ROSEN} (M.):
Remarks on the history of Fermat's last theorem 1844 to 1984,
<in Modular forms and Fermat's last theorem, Springer, New York, 1997,
pp.~505--525.>
\newcount\rosen \global\rosen=\refno

\ref{\pc WASHINGTON} (L.):
Introduction to cyclotomic fields. 
<Springer-Verlag, New York, 1997, xiv+487 pp.> 
\newcount\washington \global\washington=\refno

} 

\centerline{\bf Primary units in cyclotomic fields}
\bigskip\bigskip 
\centerline{Chandan Singh Dalawat} 
\centerline{\it Harish-Chandra Research Institute}
\centerline{\it Chhatnag Road, Jhunsi,   Allahabad 211019, India} 
\centerline{\it dalawat@gmail.com}

\bigskip\bigskip

{\pc R{\'E}SUM{\'E}}.  Nous {\'e}tudions les relations mutuelles entre
  trois notions d'unit{\'e}s $p$-primaires dans le corps cyclotomique local des
  racines $p$-i{\`e}me de~$1$ ($p$ {\'e}tant un nombre premier impair),
  sp{\'e}cialement en r{\'e}f{\'e}rence aux unit{\'e}s globales.

\bigskip

{\pc ABSTRACT}.  We investigate the interrelationships among three notions of
$p$-primary units in the local cyclotomic field of $p$-th roots of~$1$ ($p$
being an odd prime number), especially with reference to global
units.\footnote{}{MSC 2000~: 11R18, 11S10.  Keywords~: $p$-primary numbers,
  {\it nombres primaires}, {\it Prim{\"a}rzahlen}, cyclotomic fields, {\it
    corps cyclotomiques}, {\it Kreisteilungsk{\"o}rper}.}

\bigskip\bigskip

\def\primar{{\it prim{\"a}r\/}}

Let $p$ be a prime number (including $p=2$), $K$ a finite extension of $\Qp$
containing a primitive $p$-th root $\zeta$ of~$1$, $\ogoth$ the ring of
integers of $K$, and $\pgoth$ its unique maximal ideal.  Denote by
$U_n=\Ker(\ogoth^\times\to(\ogoth/\pgoth^n)^\times)$ the filtration by units
of various levels, and by $\bar U_n$ the image of $U_n$ in
$K^\times\!/K^{\times p}$.

Recall that a unit $\alpha\in\ogoth^\times$ is called $p$-{\it primary\/} if
the extension $K(\!\root p\of\alpha)$ is unramified over $K$.
It is known that $\alpha$ is $p$-primary if and only if its image in
$K^\times\!/K^{\times p}$ lies in the $\Fp$-line $\bar U_{pe_1}$
\citer\locdisc(prop.~16), where $e_1$ is the ramification index of
$K|\Q_p(\zeta)$.  This is equivalent to requiring that $\alpha$ be a $p$-th
power in $(\ogoth/\pgoth^{pe_1})^\times$ \citer\locdisc(prop.~45).

\medskip

Assume henceforth that $p$ is {\it odd}, that $K=\Qp(\zeta)$, so that $e_1=1$,
$\ogoth=\Zp[\zeta]$, and $\pgoth=\pi\ogoth$, where $\pi=1-\zeta$.  In this
case, there are two other notions of ``primary'' units, which we have called
{\it primaire\/} (\S2) and \primar\ (\S4), in order to distinguish them from
the notion of $p$-primary numbers recalled above.

The purpose of this Note is to compare these three notions, with special
reference to the global units $\Z[\zeta]^\times$.  We will show that for
$\alpha\in\Z[\zeta]^\times$, these notions are equivalent to being a $p$-th
power in $\Zp[\zeta]^\times$ (prop.~3, prop.~6), although they are
inequivalent for local units in general.  

This allows us to reconcile three different formulations of Kummer's lemma to
the effect that if the odd prime $p$ is {\it regular} --- if $p$ does not
divide the class number of $\Q(\zeta)$ --- then certain units
$u\in\Z[\zeta]^\times$ are $p$-th powers in $\Q(\zeta)^\times$.  The
difference lies in the hypotheses on~$u$~; in \citer\rosen(p.~513), $u$ is
required to be $p$-primary (at $\pgoth$)~; in \citer\bosha(p.~377), $u$ is
required to be {\it primaire\/} (\S2)~; in \citer\zahlbericht(p.~288), $u$ is
required to be \primar\ (\S4).

\bigskip

{\bf 1. $p$-primary numbers}\pointir Recall that a $1$-unit $\alpha\in U_1$ is
a $p$-th power if and only if $\alpha\in U_{p+1}$ \citer\locdisc(prop.~30)~;
$\alpha$ is $p$-primary if and only if $\alpha\in U_p$, and, finally,
$\alpha\equiv1\pmod p$ if and only if $\alpha\in U_{p-1}$.

\bigskip

{\bf 2. Nombres primaires}\pointir Traditionally, a global unit
$u\in\Z[\zeta]^\times$, or more generally an integer $u\in\Z[\zeta]$ prime to
$\pi$, is called ``primary'' if $u\equiv a\pmod p$ for some $a\in\Z$ (prime to
$p$), but the definition makes sense for all local units.  In order to
distinguish it from the notion in \S1, we will call a local unit
$\alpha\in\ogoth^\times$ {\it primaire\/} if $\alpha\equiv a\pmod p$ for some
$a\in\Z_p^\times$.  Such units form a subgroup of $\ogoth^\times$ containing
$\ogoth^{\times p}$ (lemma~1).

We show that if a global unit is {\it primaire}, then it is not only
$p$-primary in $K=\Qp(\zeta)$ but even a $p$-th power in $K^\times$.  (At the
other primes $\goth l$ of $\Q(\zeta)$, which are prime to~$p$, every global
unit $u$ is $p$-primary in the sense that adjoining $\root p\of u$ to the
local field ${\bf Q}(\zeta)_{\goth l}$ gives an unramified extension thereof,
although $u$ need not be a $p$-th power in ${\bf Q}(\zeta)_{\goth l}^\times$.)

Not every {\it primaire\/} local unit $\alpha\in\ogoth^\times$ is $p$-primary.
Indeed, we have $(\ogoth/p\ogoth)^{\times p}=\F_p^\times$ (lemma~1),
$p\ogoth=\pgoth^{p-1}$, and, in the notation of \citer\locdisc(),
$$
\hat\alpha\in(\ogoth/p\ogoth)^{\times p}
\ \ \Longleftrightarrow\ \ 
\bar\alpha\in\bar U_{p-1}
$$ 
\citer\locdisc(prop.~45), whereas $\bar U_p\neq\bar U_{p-1}$
\citer\locdisc(prop.~42). For example, $\alpha=1+p$ is $\equiv1\pmod p$ but
$\bar\alpha\notin\bar U_p$, so $1+p$ is {\it primaire\/} but not $p$-primary.

It might still be true that if a {\it global\/} unit $u\in\Z[\zeta]^\times$ is
{\it primaire\/}, then it is $p$-primary at every place of $\Q(\zeta)$, but
only the place $\pgoth|p$ really matters.  Our aim is to verify that not only
is this the case, but in fact $u\in K^{\times p}$, whether $p$ is regular or
not (prop.~3).  Let us begin with a lemma which has already been invoked.

\th LEMMA 1
\enonce
With the above notation, $(\ogoth/p\ogoth)^{\times p}=\F_p^\times$.
\endth
This is well-known, see \citer\kks(p.~130).  The inclusion
$\F_p^\times\subset(\ogoth/p\ogoth)^{\times p}$ is clear, for
$\F_p^\times=\F_p^{\times p}$.  To see the converse $(\ogoth/p\ogoth)^{\times
  p}\subset\F_p^\times$, note that the $\Fp$-space $\ogoth/p\ogoth$ admits the
basis $1,\hat\zeta,{\hat\zeta}^2,\ldots,{\hat\zeta}^{p-2}$.  Therefore, for
every element $z=\sum_{i=0}^{p-2}a_i{\hat\zeta}^i$ of $\ogoth/p\ogoth$ (with
$a_i\in\Fp$), the $p$-th power $z^p=\sum_{i=0}^{p-2}a_i$ is in $\Fp$~:
$$\eqalign{
\left(a_0+a_1\hat\zeta+\cdots+a_{p-2}{\hat\zeta}^{p-2}\right)^p
&\equiv a_0^p+a_1^p+\cdots+a_{p-2}^p\pmod{p}\cr
&\equiv a_0+a_1+\cdots+a_{p-2}\pmod{p}.\cr
}
$$

Recall that $p$ is an odd prime, that $K=\Qp(\zeta)$, that
$\ogoth=\Zp[\zeta]$, and that
$U_i=\Ker(\ogoth^\times\to(\ogoth/\pgoth^i)^\times$) for $i>0$.  The local
ingredient in Kummer's lemma amounts to $U_{p-1}\cap N\subset U_p$, where
$N=\Ker(N_{K|\Qp}:K^\times\to\Qp^\times)$.  More precisely,

\th PROPOSITION 2
\enonce
If a unit\/ $\alpha\in\ogoth^\times$ is\/ $\equiv a\pmod p$ for some\/
$a\in\Zp^\times$, and if its absolute norm\/ $N_{K|\Qp}(\alpha)$ is
$\equiv1\pmod{p\pi}$,  then\/ $\alpha$ is\/ $p$-primary.  
\endth

Notice first that we may replace $\alpha$ by $\alpha^{p-1}$~: adjoining
$\root p\of\beta$ or $\root p\of{\beta^{p-1}}$ gives the same extension of
$K$, for any $\beta\in K^\times$.  We may thus assume that $\alpha\equiv1\pmod
p$, and write $\alpha=1+\gamma p$ for some $\gamma\in\ogoth$, and, as $\ogoth$
and $\Zp$ have the same residue field, $\gamma=c+\delta\pi$ for some $c\in\Zp$
and $\delta\in\ogoth$, so that $\alpha=1+cp+\delta p\pi$.  Now, for every
$\sigma\in\Gal(K|\Qp)$, we have
$\sigma(\alpha)=1+cp+\sigma(\delta)p\sigma(\pi)$, so that
$\sigma(\alpha)\equiv1+cp\pmod{p\pi}$, for $\sigma(\pi)$ is also a uniformiser
of $K$.  Taking the product over all $\sigma$, we get
$$
1\equiv N_{K|\Qp}(\alpha)\equiv(1+cp)^{p-1}\equiv1-cp\pmod{p\pi}.
$$
This implies that $cp\equiv0\pmod{p\pi}$, and hence
$\alpha\equiv1\pmod{p\pi}$, showing that $\alpha$ is $p$-primary
\citer\locdisc(prop.~16).  This proof is adapted from
\citer\washington(p.~80).

\smallskip

\th PROPOSITION 3
\enonce
If a global unit\/ $u\in\Z[\zeta]^\times$ is {\it primaire\/}, then it is a
$p$-th power in\/ $\ogoth^\times$.
\endth
[If $\bar E$ is the image of the global units $\Z[\zeta]^\times$ in the group
$\overline{\ogoth^\times}=\ogoth^\times\!/\ogoth^{\times p}$ of local units
modulo $p$-th powers, and if $(\bar U_n)_{n>0}$ denotes the filtration on the
latter group, so that $\overline{\ogoth^\times}=\bar U_1$ and $\bar
U_{p+1}=\{1\}$, then $\bar E\cap\bar U_{p-1}=\{1\}$, although $\bar U_{p-1}$
is $2$-dimensional over $\Fp$.]

We first need to recall a few facts about global units.  Every
$u\in\Z[\zeta]^\times$ is (uniquely) of the form $u=\xi w$ for some $p$-th
root $\xi$ of~$1$ and some $w\in\Z[\zeta+\zeta^{-1}]^\times$ ``totally real''
(\citer\washington(p.~3), \citer\kks(p.~129)).  If moreover $u\equiv
a\pmod{\pi^2}$ for some $a\in\Z$, then $u=w$, for $\xi\in U_1$ but $\xi\notin
U_2$, unless $\xi=1$ \citer\washington(p.~79).  Hence
$N_{\Q(\zeta)|\Q}(u)=N_{\Q(\zeta+\zeta^{-1})|\Q}(u)^2=1$, for the norm of a
unit is a unit and $\Z^{\times2}=\{1\}$.  Also,
$N_{K|\Qp}(u)=N_{\Q(\zeta)|\Q}(u)=1$.

Now suppose that $u$ is {\it primaire}~; in particular, $u\equiv
a\pmod{\pi^2}$ for some $a\in\Z$.  The above discussion implies that
$N_{K|\Qp}(u)=1$, and prop.~2 then implies that $u$ is $p$-primary.  

But the above discussion also implies that $u$ is in
$K^+=\Qp(\zeta+\zeta^{-1})$.  Up to multiplying $u$ by a $p$-th power (such as
the multiplicative representative $\langle a^{-1}\rangle\in\ogoth^\times$,
where $a\in\F_p^\times$ is the image of $u$), or replacing $u$ by $u^{p-1}$,
we may assume that $u\in U_1$.  As for any $p$-primary $1$-unit of $K$, we
have $u\in U_p$.  But if a unit of $K^+$ (such as $u$) is in $U_n$ for some
{\it odd\/}~$n$ (such as $n=p$), then it is in $U_{n+1}$, because the
ramification index of $K|K^+$ is~$2$.  But $U_{p+1}=U_2^p$, so it follows that
$u\in K^{\times p}$ \citer\locdisc(prop.~30).

%

{\it Remark}\pointir Prop.~3 allows us to prove Kummer's lemma in the variant
\citer\bosha(p.~377) along the same lines as in \citer\rosen(p.~513), thereby
avoiding the $p$-adic logarithm or any extraneous global considerations.
Namely, if $u\in\Z[\zeta]^\times$ is {\it primaire}, then it is $p$-primary,
even a local $p$-th power (prop.~3), and hence the extension obtained by
adjoining $\root p\of u$ to $\Q(\zeta)$ is cyclic (of degree~$1$ or~$p$) and
unramified everywhere.  But if $p$ is {\it regular}, $\Q(\zeta)$ has no
everywhere-unramified cyclic degree-$p$ extension, by class field theory (or
as a consequence of Hilbert's {\it Satz\/}~94, as in \citer\rosen(p.~523)).
Hence $u\in\Z[\zeta]^{\times p}$.

\medbreak

{\bf 3. A general observation}\pointir The local argument of prop.~2 can be
generalised so as to bring out its essential features.  Let $F$ be {\it any\/}
finite extension of $\Qp$, allow $p$ to be~$2$, and let $L|F$ be a totally but
tamely ramified extension of degree~$e$ (prime to~$p$).  Let $\pi_F$ and
$\pi_L$ be uniformisers of $F$, $L$ respectively.  If a unit
$\alpha\in\ogoth_L^\times$ is $\equiv a\pmod{\pi_F^r}$ for some
$a\in\ogoth_F^\times$ and some $r>0$, then clearly its relative norm
$N_{L|F}(\alpha)$ is $\equiv a^e\pmod{\pi_F^r}$.  But if we demand that
$N_{L|F}(\alpha)\equiv a^e\pmod{\pi_F^r\pi_L^{\phantom{r}}}$, then it follows
that $\alpha\equiv a\pmod{\pi_F^r\pi_L^{\phantom{r}}}$, as the following
prop.\ shows.

\th PROPOSITION 4
\enonce
Suppose that\/ $L|F$ is totally ramified of degree\/~$e$ prime to\/~$p$, and
let\/ $\alpha\in\ogoth_L^\times$.  If\/ $\alpha\equiv a\pmod{\pi_F^r}$
and\/ $N_{L|F}(\alpha)\equiv a^e\pmod{\pi_F^r\pi_L^{\phantom{r}}}$ for some\/
$a\in\ogoth_F^\times$ and some\/ $r>0$, then\/ $\alpha\equiv
a\pmod{\pi_F^r\pi_L^{\phantom{r}}}$. 
\endth
Write $\alpha=a+b\pi_F^r+\gamma\pi_F^r\pi_L^{\phantom{r}}$, where we assume
that $b\in\ogoth_F$ because $L|F$ is totally ramified, and
$\gamma\in\ogoth_L$.  For every $F$-conjugate $\sigma(\alpha)$ of $\alpha$, we
have $\sigma(\alpha)=a+b\pi_F^r+\sigma(\gamma)\pi_F^r\sigma(\pi_L)$, and,
taking the product over the $e$ $F$-embeddings $\sigma$ of $L$ (in some fixed
algebraic closure of $F$), we get
$$
a^e\equiv N_{L|F}(\alpha)\equiv 
(a+b\pi_F^r)^e\equiv 
a^e+ea^{e-1}b\pi_F^r\pmod{\pi_F^r\pi_L^{\phantom{r}}}.
$$
Therefore, working $\!\!\pmod{\pi_F^r\pi_L^{\phantom{r}}}$, we have
$ea^{e-1}b\pi_F^r\equiv0$, so  $b\pi_F^r\equiv0$ (as $e$ and $a$ are units in
$L$) and hence $\alpha\equiv a$. 

Let $U_m$ (resp.~$N_m$) be the group of $\alpha\in\ogoth_L^\times$ such that
$\alpha$ (resp.~$N_{L|F}(\alpha))$ is $\equiv1\pmod{\pi_L^m}$.  Taking $a=1$
in prop.~4, we get

\th COROLLARY 5
\enonce
When\/ $L|F$ is totally tamely ramified of degree\/~$e$, we have\/
$U_{re}\cap N_{re+1}=U_{re+1}$ for every\/ $r>0$.
\endth

Prop.~2 was essentially the case $F=\Qp$ ($p$ odd), $L=\Qp(\zeta)$, $r=1$.
So, from our local perspective, the basic point in Kummer's lemma in the
formulation \citer\bosha(p.~377) is that this $L|F$ is totally (but tamely)
ramified of degree~$p-1$, and that if $u\in\Z[\zeta]^\times$ is {\it
primaire\/}, then its norm is~$1$.  Therefore $u$ is $p$-primary at $\pi$ (and
also at every other place of $\Q(\zeta)$).

\medskip

{\bf 4. Prim{\"a}rzahlen}\pointir In Hilbert's {\it Zahlbericht\/}, there is a
third, closely allied, notion.  He defines it only for prime-to-$\pi$ integers
in $\Z[\zeta]$, but it makes sense for all local units
$\alpha\in\ogoth^\times$, where $\ogoth=\Zp[\zeta]$ (and $p$ is an {\it odd\/}
prime).  We say that $\alpha$ is \primar\ if
$$
\alpha\equiv a\pmod{\pi^2},\quad 
N_{K|K^+}(\alpha)\equiv b\pmod p,
\qquad (a,b\in\Z_p^\times),
$$
where $K$, $K^+$ are the completions of $\Q(\zeta)$, $\Q(\zeta+\zeta^{-1})$ at 
the unique place above~$p$.  Such units form a subgroup of $\ogoth^\times$~; 
the name has been chosen to distinguish them from $p$-primary (\S1) or
{\it primaire\/} (\S2) units.

It is clear that if a local unit $\alpha\in\ogoth^\times$ is {\it primaire},
then it is \primar, for $\alpha\equiv a\pmod p$ implies $\alpha\equiv
a\pmod{\pi^2}$ and $N_{K|K^+}(\alpha)\equiv a^2\pmod p$.  In particular, every
$p$-th power $\alpha\in\ogoth^{\times p}$ is \primar\ (lemma~1).

The converse is of course true for $p=3$, but false for $p\neq3$.  Indeed, let
$\varpi$ be a $(p-1)$-th root of $-p$ in $K$ \citer\locdisc(prop.~24), so that
$\varpi$ is a uniformiser of $K$ and $\sigma_{-1}(\varpi)=-\varpi$, where
$\sigma_{-1}$ is the generator of $\Gal(K|K^+)$.  It is clear that
$\alpha=1+\varpi^{p-2}$ is \primar\ but not {\it primaire}, if $p>3$.

If a global unit $u\in\Z[\zeta]^\times$ is \primar, and if the prime $p$ is
{\it regular}, then $u\in\Z[\zeta]^{\times p}$ \citer\zahlbericht(p.~288).
One might therefore suspect that, in general, a \primar\ global unit is
$p$-primary (at the prime $\pi$), even if $p$ is irregular.  We show that in
fact a \primar\ global unit is always a $p$-th power in $K^\times$.

\th PROPOSITION 6
\enonce
If a global unit\/ $u\in\Z[\zeta]^\times$ is \primar, then\/
$u\in\ogoth^{\times p}$.
\endth
[Denoting by $\bar P$ the image in
$\overline{\ogoth^\times}=\ogoth^\times\!/\ogoth^{\times p}$ of the group of
\primar\ local units, and by $\bar E$ the image of all global units
$\Z[\zeta]^\times$, we have $\bar P\cap\bar E=\{1\}$.]

Suppose that $u$ is \primar~; it is enough to show that $u'=u^{p-1}$ is a
$p$-th power in $K^\times$.  Since $u'\equiv 1\pmod{\pi^2}$, we have
$u'\in\Z[\zeta+\zeta^{-1}]^\times$, as in the proof of prop.~3, and hence
$u'\in K^+$.  But then $N_{K|K^+}(u')=u^{\prime2}$, and, as $N_{K|K^+}(u')\in
U_{p-1}$ by hypothesis, we have $u'\in U_{p-1}$, for $(\ )^{1/2}$ is an
automorphism of the $\Zp$-module $U_{p-1}$.

We have shown that $u'$ is {\it primaire}.  As it is also a global unit,
prop.~3 implies that $u'\in\ogoth^{\times p}$.  Hence $u\in\ogoth^{\times p}$.

(The local ingredient in the above proof says more generally that if
$\alpha\in U_{p-1}\cap K^+$ and if $N_{K|\Qp}(\alpha)\equiv1\pmod{p\pi}$, then
$\alpha\in U_{p+1}\subset\ogoth^{\times p}$.  Indeed, $\alpha\in U_p$ by
prop.~2, and $U_p\cap K^+\subset U_{p+1}$, because the ramification index of
$K|K^+$ is~$2$ and $p$ is odd.)

\bigskip

{\bf 5. Summary}\pointir Let us summarise.  For local units
$\alpha\in\ogoth^\times$, we have
$$
\alpha\in\ogoth^{\times p}\ \Longrightarrow\ 
\alpha\hbox{ is $p$-primary}\ \Longrightarrow\ 
\alpha\hbox{ is {\it primaire}}\ \Longrightarrow\ 
\alpha\hbox{ is \primar},
$$
and all these implications are strict, except that the last one is an
equivalence when $p=3$.  But for global units $\alpha\in\Z[\zeta]^\times$, all
three implications are actually equivalences (prop.~3, prop.~6).

Let $\bar P$ be the image in
$\overline{\ogoth^\times}=\ogoth^\times\!/\ogoth^{\times p}$ of the group of
\primar\ local units~; we have $\bar U_{p-1}\subset\bar P\subset\bar U_2$.
The above implications can be rewritten as
$$
\{1\}\subset\bar U_p\subset\bar U_{p-1}\subset\bar P
$$ 
where all three inclusions are strict, except for the last one, which is an
equality when $p=3$.  Finally, $\bar E\subset\overline{\ogoth^\times}$ being
the image of the global units $\Z[\zeta]^\times$, we have $\bar E\cap\bar
P=\{1\}$ (prop.~6).  In short, although the four notions are distinct locally
(at $\pi$), they are equivalent globally.

(Notice that the $\Fp$-dimension of $\bar P/\bar U_{p-1}$ grows linearly
with~$p$, and equals the number of odd $a\in[3,p-2]$ such that $2a\ge p-1$.
Indeed, denoting the set of such $a$ by $I$, a basis is provided by the images
of $1+\varpi^a$ ($a\in I$), where $\varpi^{p-1}=-p$.)

\medbreak

{\bf 6. A suggestion}\pointir Starting with Kummer, it is proved at many
places that if $u\in\Z[\zeta]^\times$ is $p$-primary (\S1) or {\it primaire\/}
(\S2) or \primar\ (\S4), and if $p$ is regular, then $u\in\Z[\zeta]^{\times
p}$.  We are advocating that this be done in two steps.  The first step, which
we have carried out here, is essentially local and says that if
$u\in\Z[\zeta]^\times$ is $p$-primary or {\it primaire\/} or merely \primar,
then $u\in\Z_p[\zeta]^{\times p}$~; it is valid for all odd primes $p$,
regular or not.  The second step, which is global, says of course that if
moreover $p$ is regular, then $u\in\Z[\zeta]^{\times p}$.

\bigskip {\it Acknowledgement.}  The author is grateful to Professor K N
Ganesh, director of the Indian Institute of Science Education and Research,
Pune, for having gone to great lengths to make a truly enlightening remark on
24/9/2009.

\bigbreak
\unvbox\bibbox 

\bye